\documentclass[12pt]{article}
\usepackage{amssymb}
\usepackage{amsmath}
\usepackage{array}
\usepackage[T2A]{fontenc}
\usepackage[cp1251]{inputenc}
\usepackage[russian]{babel}
\usepackage[pdftex]{graphics}

\DeclareMathOperator{\RRe}{Re} \DeclareMathOperator{\IIm}{Im}

\DeclareMathOperator{\vep}{\varepsilon}

\DeclareMathOperator{\vf}{\varphi}

\newcommand{\Z}{\scriptstyle}
\newcommand{\D}{\displaystyle}

\multlinegap=0em \DeclareFontFamily{T1}{msb}{}
\DeclareFontShape{T1}{msb}{m}{ol}{<5> <6> <7> <8> <9> gen * msbm
<10> <10.95> <12> <14.4> <17.28> <20.74> <24.88> msbm10}{}
\DeclareSymbolFont{AMSb}{T1}{msb}{m}{ol} \multlinegap=0em

\renewcommand{\S}{\mathhexbox278}
\renewcommand{\le}{\operatorname{\leqslant}}
\renewcommand{\ge}{\operatorname{\geqslant}}

\begin{document}

\begin{flushleft}
11M06
\end{flushleft}

\begin{center}
{\rmfamily\bfseries\normalsize M.A.~Korolev\footnote{The author is
supported by RFBR (grant no. 12-01-33080).}}
\end{center}

\begin{center}
{\rmfamily\bfseries\normalsize On the large values of the Riemann
zeta\,-function on the critical line - II}
\end{center}

\vspace{0.5cm}

\fontsize{11}{12pt}\selectfont

\textbf{Annotation.} We prove some new bounds for the maximum of
$\bigl|\zeta(0.5+it)\bigr|$ on the segments $T\le t\le T+H$ with
$H\,\ge\,(\ln\ln\ln{T})^{1+\vep}$. All the theorems are based on the
Riemann hypothesis.

\vspace{0.2cm}

\fontsize{12}{15pt}\selectfont

\vspace{0.5cm}

\textbf{Introduction}

We continue the investigation of the lower bound estimates for the
maximum of modulus of the Riemann zeta\,-function $\zeta(s)$ on the
short segments of the critical line $\RRe s =0.5$.

The theorem of R.~Balasubramanian \cite{Balasubramanian_1986} states
that the function
\[
F(T;H)\,=\,\max_{|t-T|\le H}\bigl|\zeta(0.5+it)\bigr|
\]
satisfies the inequality
\begin{equation}\label{lab01}
F(T;H)\,\gg\,\exp{\biggl(\frac{3}{4}\sqrt{\frac{\ln{H}}{\ln\ln{H}}}\biggr)}
\end{equation}
for $\ln\ln{T}\ll H\le 0.1T$. It is supposed that this bound is
close to the best possible (at least, for $H\asymp T$; see
\cite{Farmer_Gonek_Hughes_2007})). In the case of ``very small''
$H$, $0<H\ll \ln\ln{T}$, there is a series of lower bound estimates
for $F(T;H)$, but all of them differ essentially from (\ref{lab01}),
because their right hand side decreases when $T$ grows (see
\cite{Karatsuba_2001a}-\cite{Changa_2005}).

In particular, it was proved in \cite{Karatsuba_2006} that
\begin{equation}\label{lab02}
F(T;H)\,\ge\,\frac{1}{16}\,\exp\biggl\{-\,\frac{5\ln{T}}{6(\pi/\alpha\,-\,1)(\ch{(\alpha
H)-1})}\biggr\}
\end{equation}
for any fixed $\alpha$, $1\le \alpha <\pi$, $2\le\alpha H\le
\ln\ln{T}-c_{1}$, where $c_{1}>0$ is some absolute constant. Given
$\vep>0$, it follows from (\ref{lab02}) that for any $T\ge
T_{0}(\vep)>0$ and for $H\ge \pi^{-1}(1+\vep)\ln\ln{T}-c_{1}$, the
function $F(T;H)$ is bounded from below by some constant:
\[
F(T;H)\,>\,c_{2}\,=\,\frac{1}{16}\,\exp{\bigl(-\,1.7\vep^{-1}e^{c_{1}}\bigr)}>0.
\]
In \cite{Karatsuba_2006}, A.A.~Karatsuba posed the problem of
proving $F(T;H)\ge 1$ for the values of $H$ essentially smaller than
$\ln\ln{T}$, namely, for $H\ge \ln\ln\ln{T}$. The conditional
solution of this problem was obtained in \cite{Korolev_2014}.
Namely, it was proved that for an arbitrary large but fixed constant
$A>1$ there exist (non-effective) constants $c_{0},T_{0}$ such that
\begin{equation}\label{lab03}
F(T;H)\ge A\quad \text{for any} \quad T\ge T_{0}\quad
\text{and}\quad H\ge \pi^{-1}\ln\ln\ln{T}+c_{0}.
\end{equation}

The comparison of (\ref{lab01}) and (\ref{lab03}) leads us to the
following questions:

\vspace{0.2cm}

1) for what size of $H$, $H\ll \ln\ln{T}$, the inequality
$F(T;H)\,\gg\,\exp{\bigl((\ln{H})^{0.5-\vep}\bigr)}$ holds?

\vspace{0.2cm}

2) for what size of $H$, $H\gg \ln\ln\ln{T}$, the inequality
$F(T;H)\,\gg\,f(H)$ holds for some unbounded function $f(u)$?

\vspace{0.2cm}

The main goal of the present paper is to prove the following
assertions based on the Riemann hypothesis (RH).

\vspace{0.2cm}

\textsc{Theorem 1.} \emph{Suppose that RH is true, and let $m\ge 1$
is any fixed integer. Then}
\[
F(T;H)\,\ge\;\exp{\biggl(\frac{0.05\sqrt{\ln{H}}}{(2m\ln\ln{H})^{m}}\biggr)}
\]
\emph{for} $T\ge T_{0}$ \emph{and} $(\ln\ln{T})^{\frac{\Z 1}{\Z
2m\mathstrut}}\,\le\,H\,\le\,\ln\ln{T}$.

\vspace{0.2cm}

\textsc{Theorem 2.} \emph{Suppose that RH is true, and let
$0<\vep<0.1$ be any fixed number. Then}
\[
F(T;H)\,\ge\,\exp{\left(\sqrt{\ln{H}}\,e^{-c\,(\ln\ln{H})^{\Z
1-0.5\vep}}\right)}
\]
\emph{for any $T\ge T_{1}(\vep)$, $H\ge\,(\ln\ln\ln{T})^{2+\vep}$
and for some constant} $c = c(\vep)>0$.

\vspace{0.2cm}

\textsc{Theorem 3.} \emph{Suppose that RH is true, and let
$0<\vep<0.1$ be any fixed number. Then}
\[
F(T;H)\,\ge\,\exp{\left((\ln{H})^{\gamma\,-\,\vep}\right)}
\]
\emph{for any $T\ge T_{1}(\vep)$, $H\ge\,(\ln\ln\ln{T})^{2}$,
\[
\gamma\,=\,\frac{1}{2+(\pi\varrho)^{-1}}\,=\,0.46862145\ldots,
\]
where $\varrho = 2.37689234\ldots$ stands for the least positive
root of the function}
\[
h(\lambda)\,=\,\int_{0}^{+\infty}e^{-(\ch{\sqrt{u}}+\cos{\sqrt{u}}\,)}\cos{(\lambda
u)}\,du.
\]

\vspace{0.2cm}

\textsc{Theorem 4.} \emph{Suppose that RH is true, and let
$0<\vep<0.1$ be any fixed number. Then}
\[
F(T;H)\,\ge\,\exp{\left(0.5\,e^{(\ln\ln{H})^{\Z 0.5\vep}}\right)}
\]
\emph{for any} $T\ge T_{1}(\vep)$, $H\ge\,(\ln\ln\ln{T})^{1+\vep}$.

\vspace{0.2cm}

The proof of all the above assertions is based on the general
Theorem A. Its particular cases are used in
\cite{Selberg_1946}-\cite{Korolev_2013}. At the same time, the proof
of Theorem A is based on the convolution formula (lemma 1 of present
paper) going back to A.~Selberg (see \cite{Selberg_1946} and
\cite{Tsang_1986}) and on the lemma of prof. K.-M.~Tsang (see lemma
2 below). The original parts of paper are the upper bound estimates
for the rate of decreasing for Fourier transforms of some rapidly
decreasing functions (lemma 4). The idea of varying of the
function $f(u)$ in convolution formula for minimization of $H$
belongs to R.N.~Boyarinov \cite{Boyarinov_2011a}, \cite{Boyarinov_2011b}.

\vspace{0.5cm}

\textbf{\S 1. Auxilliary assertions}

\vspace{0.2cm}

In this section, we give some auxilliary assertions needed for the
proof of Theorem A.

\vspace{0.2cm}

\textsc{Lemma 1.} \emph{Suppose that the function $f(z)$ is
analytical in the strip $|\!\IIm z|\leqslant 0.5+\alpha$, where it
satisfies the inequality $|f(z)|\leqslant c(|z|+1)^{-(1+\beta)}$
with some positive $\alpha, \beta$ and $c$. Then the identity}
 \begin{multline*}
 \int_{-\infty}^{+\infty}f(u)\ln{\zeta\bigl(0.5+i(t+u)\bigr)}\,du\,=\,
 \sum\limits_{n = 2}^{+\infty}\frac{\Lambda_{1}(n)}{\sqrt{n}}\,n^{-it}\widehat{f}(\ln{n})\,+\\
 +\,2\pi\biggl(\;\sum\limits_{\beta>0.5}\int_{0}^{\beta-0.5}f(\gamma-t-iv)\,dv\,-\,\int_{0}^{0.5}f(-t-iv)\,dv\biggr),
 \end{multline*}
\emph{holds for any $t$, where $\varrho=\beta+i\gamma$ in the last
sum runs through all complex zeros of $\zeta(s)$ to the right from
the critical line, $\Lambda_{1}(n)=\Lambda(n)/\ln{n}$.}

\vspace{0.2cm}

This assertion goes back to A.~Selberg (see for example \cite[Lemma
16]{Selberg_1946}). In \cite[Ch. II, \S 11]{Tsang_1984}, \cite[Ch.
II, \S 2]{Karatsuba_Korolev_2005}, \cite{Tsang_1986}, there are some
variants of this lemma, where $f(z)$ satisfies slightly different
conditions. These proofs can be easily adopted to the case under
considering.

\vspace{0.2cm}

\textsc{Lemma 2.} \emph{Let $H>0, M>0$, $k\ge 1$, and suppose that the
real function $W(t)$ satisfies to the following conditions:
\[
\int_{T}^{T+H}W^{2k}(t)\,dt\,>\,HM^{2k},\qquad
\biggl|\int_{T}^{T+H}W^{2k+1}(t)\,dt\biggr|\,\le\,0.5HM^{2k+1}.
\]
Then}
\[
\max_{T\le t\le T+H}\bigl(\pm W(t)\bigr)\,\ge\,0.5M.
\]

\vspace{0.2cm}

This is a small modification of lemma 11.3 from \cite[Ch. II, \S
11]{Tsang_1984} (see also lemma 4 from \cite{Tsang_1986}).

\vspace{0.2cm}

\textsc{Lemma 3.} \emph{Let $\varpi_{1} = 2, \varpi_{2}=3,
\varpi_{3} = 5, \ldots,\ldots$ are all the primes indexed in
ascending order. Then $\varpi_{n}<n(\ln{n}+\ln\ln{n})$ for any $n\ge
6$. Next,
\[
\sum\limits_{\varpi_{n}\le
x}\frac{1}{\varpi_{n}}\,=\,\ln\ln{x}\,+\,\mathfrak{m}\,+\,\frac{\theta}{\ln^{2}x},
\]
where $\mathfrak{m} = 0.261 497 \ldots$ is Mertens' constant and
$-0.5<\theta<1$ for any $x>1$.}

\vspace{0.2cm}

These assertions follows from theorems 3, 5 and 6 of
\cite{Rosser_Schoenfeld_1962}.

\vspace{0.5cm}

\textbf{\S 2. General theorem}

\vspace{0.2cm}

This section is devoted entirely to the  proof of one general
assertion, which implies all the theorems 1-4.

\pagebreak

\textsc{Theorem A.} \emph{Suppose RH is true, and let the function
$\Phi(u)$ satisfies the following conditions:}

1) \emph{$\Phi(u)\ge 0$ for real $u$ and $f(z) = \Phi(\tau z)$ is
analytic in the strip $|\!\IIm z|\le 0.5+\delta$ for any $\tau>0$
and satisfies the inequality $|f(z)|\ll (1+|z|)^{-(1+\beta)}$ for
some positive $\beta$ and $\delta$} (\emph{both $\beta$ and $\delta$
may depend on $\tau$});

2) \emph{$|\Phi(u)|\le e^{-G(|u|)}$ for any real $u$, $|u|\ge
u_{0}$, where the functions $G(u)$, $G'(u)$ are positive and
unboundedly increasing and such that the functions $g\,'(v),
g\,'(v)\ln{g(v)}$ are positive and decreasing for $v\ge v_{0}>0$}
(\emph{here $g(v)$ stands for the inverse function to $G(u)$});

3) \emph{$\widehat{\Phi}(\lambda)$ is real for real $\lambda$, strictly
positive and monotonically decreasing on $[0,\alpha]$ for some
$\alpha > 0$; moreover,}
\begin{equation}\label{lab04}
|\widehat{\Phi}(\lambda)|\,\le\,e^{-\,|\lambda|F(|\lambda|)}
\end{equation}
\emph{for some increasing function $F(u)$ and for any real
$\lambda$, $|\lambda|\ge \lambda_{0}$};

4) \emph{the function $\vf(v)$, which is inverse to $F(u)$, is
increasing for $v\ge v_{0}$ and satisfies the inequalities}
\begin{equation}\label{lab05}
\ln{v}\,\le\,\ln{\vf(v)}\,\le\,e^{0.5\alpha v};
\end{equation}

\emph{Suppose also that $\tau_{0}$ is a root of the transcendental
equation}
\begin{equation}\label{lab06}
\alpha\tau_{0}\,+\,\ln\vf\biggl(\frac{\tau_{0}}{2}+1\biggr)\,=\,\ln\ln{H},
\end{equation}
\emph{which is unique when $H$ is sufficiently large. Finally, let
$T\ge T_{0}(\Phi;\alpha)>0$ and $H\tau_{0}\ge$ $g(\ln\ln{T})$. Then}
\begin{equation}\label{lab07}
\max_{T-H\le t\le T+2H}\ln{\bigl|\zeta(0.5+it)\bigr|}\,\ge\,\mu^{*},
\end{equation}
\emph{where}
\[
\mu^{*}\,=\,\frac{1}{10\alpha}\,\frac{\widehat{\Phi}(\alpha)}{\widehat{\Phi}(0)}\,\sqrt{\frac{\ln{\kappa}}{\kappa}}\,
\frac{e^{0.5\alpha\tau_{0}}}{\tau_{0}},
\quad \kappa = \max{(61,4\alpha^{-1})}.
\]
\emph{If, in addition, the function $\vf(v)$ satisfies the
condition}
\begin{equation}\label{lab08}
\ln{\vf(v)}\,\le\,0.5v,
\end{equation}
\emph{and $\tau_{1}$ denotes the root of the equation}
\begin{equation}\label{lab09}
\frac{\alpha\tau_{1}}{2}\,+\,\ln\vf\biggl(\frac{\tau_{1}}{2}+1\biggr)\,=\,\ln\ln{H},
\end{equation}
\emph{then the inequality}
\begin{equation}\label{lab10}
\max_{T-H\le t\le T+2H}\ln{\bigl|\zeta(0.5+it)\bigr|}\,\ge\,\mu^{**}
\end{equation}
\emph{holds for $T_{1}(\Phi;\alpha)>0$, $H\tau_{1}\ge g(\ln\ln{T})$
with}
\[
\mu^{**}\,=\,\frac{1}{6\sqrt{\alpha\kappa}}\,\frac{\widehat{\Phi}(\alpha)}{\widehat{\Phi}(0)}\,\frac{e^{0.25\alpha\tau_{1}}}{\sqrt{\tau_{1}}},
\quad \kappa = \max{(0.5,4\alpha^{-1})}.
\]

\textsc{Corollary.} \emph{Suppose that all the conditions are
satisfied and let $\varrho > 0$ be the least positive root of
$\widehat{\Phi}(\lambda)$. If $\widehat{\Phi}(\lambda)$ decreases on
$[0,\varrho]$, then the inequality} (\ref{lab07}) \emph{holds with}
\[
\mu^{*}\,=\,\frac{1}{5e\varrho}\sqrt{\frac{\ln{\kappa}}{\kappa}}\,
\frac{|\widehat{\Phi}'(\varrho)|}{\widehat{\Phi}(0)}\,\frac{e^{\,0.5\varrho\tau_{0}}}{\tau_{0}^{2}},\quad
\kappa = \max{(32,5\varrho^{-1})},
\]
\emph{where $\tau_{0}$ denotes the root of} (\ref{lab06})
\emph{corresponding to $\alpha = \varrho$. In, in addition, the
condition} (\ref{lab08}) \emph{holds true, then the inequality}
(\ref{lab07}) \emph{is true for}
\[
\mu^{**}\,=\,\frac{1}{5\sqrt{e\kappa}}\,
\frac{|\widehat{\Phi}'(\varrho)|}{\widehat{\Phi}(0)}\,\frac{e^{\,0.25\varrho\tau_{1}}}{\tau_{1}^{1.5}},\quad
\kappa = \max{(4,0.5\varrho)},
\]
\emph{where $\tau_{1}$ denotes the root of} (\ref{lab09})
\emph{corresponding to} $\alpha = \varrho$.

\vspace{0.2cm}

\textsc{Proof.} Let $\tau_{0}$ be a root of (\ref{lab06}), and suppose that $H\tau_{0}\ge g(\ln\ln{T})$, $T\ge T_{0}(\Phi;\alpha)$.
By lemma 1, $I(t) = A(t) - B(t)$, where
\begin{multline}\label{lab11}
I(t)\,=\,\int_{-\infty}^{+\infty}\Phi(\tau_{0}u)\ln{\bigl|\zeta(0.5+i(t+u))\bigr|}du,\\
A(t)\,=\,\frac{1}{\tau_{0}}\sum\limits_{n = 2}^{+\infty}\frac{\Lambda_{1}(n)}{\sqrt{n}}\,\widehat{\Phi}\biggl(\frac{\ln{n}}{\tau_{0}}\biggr)\cos{(t\ln{n})},\\
B(t)\,=\,2\pi\int_{0}^{0.5}\RRe\Phi\bigl(-(t+iu)\tau_{0}\bigr)\,du.
\end{multline}
Transforming $I(t)$, we get
\[
I(t)\,=\,\biggl(\;\int_{-H}^{H}\,+\,\int_{H}^{+\infty}\,+\,\int_{-\infty}^{-H}\,\biggr)\ldots du\,=\,I_{0}(t)+I_{1}(t)+I_{2}(t).
\]
Estimating $I_{1}, I_{2}$ from above, we note that if $\bigl|\zeta(0.5+i(t+u))\bigr|<1$ for every $u$, $u_{1}\le u\le u_{2}$, then the integral over $(u_{1},u_{2})$ is negative.
Hence, it is sufficient to estimate the integrals over the set of $u$ such that $\bigl|\zeta(0.5+i(t+u))\bigr|\ge 1$. Thus, the trivial bound
\[
\bigl|\zeta(0.5+i(t+u))\bigr|\le |t+u|+3
\]
yields
\[
I_{1}(t)\,=\,\int_{H}^{+\infty}\Phi(\tau_{0}u)\ln{(t+u+3)}du\,\le\,\ln{(2t+3)}\int_{H}^{t}\Phi(\tau_{0}u)\,du
+ 2\int_{t}^{+\infty}\Phi(\tau_{0}u)(\ln{u})\,du.
\]
Standing $j_{1}, j_{2}$ for the last integrals, we get
\begin{multline*}
j_{1}\,=\,\frac{1}{\tau_{0}}\int_{\tau_{0}H}^{+\infty}\Phi(v)dv\,\le\,\frac{1}{\tau_{0}}\int_{\tau_{0}H}^{+\infty}e^{-G(v)}dv\,=\,
\frac{1}{\tau_{0}}\int_{G(\tau_{0}H)}^{+\infty}e^{-\,u}g'(u)du\,\le\\
\le\,\frac{g'(G(\tau_{0}H))}{\tau_{0}}\int_{G(\tau_{0}H)}^{+\infty}e^{-\,u}du\,=\,\frac{g'(G(\tau_{0}H))}{\tau_{0}}\,e^{-G(\tau_{0}H)}.
\end{multline*}
Since $g'(G(u))G'(u) = 1$ then
\[
j_{1}\,\le\,\frac{e^{-\,G(\tau_{0}H)}}{\tau_{0}G'(\tau_{0}H)}.
\]
Similarly we have
\begin{multline*}
j_{2}\,\le\,\frac{1}{\tau_{0}}\int_{\tau_{0}t}^{+\infty}\Phi(v)\biggl(\ln\frac{v}{\tau_{0}}\biggr)dv\,\le\,\frac{1}{\tau_{0}}\int_{\tau_{0}t}^{+\infty}e^{-G(v)}
\biggl(\ln\frac{v}{\tau_{0}}\biggr)dv\,=\\
=\,\frac{1}{\tau_{0}}\int_{G(\tau_{0}t)}^{+\infty}e^{-\,u}g'(u)\,\ln{\biggl(\frac{g(u)}{\tau_{0}}\biggr)}du\,\le\,
\frac{g'(G(\tau_{0}t))}{\tau_{0}}\,\ln{\biggl(\frac{1}{\tau_{0}}\,g(G(\tau_{0}t))\biggr)}e^{-\,G(\tau_{0}t)}\,=\\
=\,\frac{e^{-\,G(\tau_{0}t)}}{\tau_{0}}\,\frac{\ln{t}}{G'(\tau_{0}t)}.
\end{multline*}
Hence,
\[
I_{1}\,\le\,\frac{2(\ln{t})e^{-\,G(\tau_{0}H)}}{\tau_{0}G'(\tau_{0}H)}.
\]
The same bound holds for $I_{2}$. Thus,
\begin{equation}\label{lab12}
I(t)\,\le\,I_{0}\,+\,\frac{4(\ln{t})}{\tau_{0}}\,\frac{e^{-\,G(\tau_{0}H)}}{G'(\tau_{0}H)}.
\end{equation}
Further we have
\begin{equation}\label{lab13}
|B(t)|\,\le\,2\pi\int_{0}^{0.5}\frac{c\,du}{(1+|t+iu|)^{1+\beta}}\,\le\,\frac{\pi
c}{t}.
\end{equation}
We split the sum $A(t)$ to the parts $A_{1}, A_{2}$ and $A_{3}$ according to the conditions $p\le X$, $n=p^{k}\le X, k\ge 2$ ($p$ is prime) and $n>X$, where
\[
X\,=\,\exp{\biggl(\tau_{0}\vf\biggl(\frac{\tau_{0}}{2}+1\biggr)\biggr)}.
\]
First we have
\begin{multline}\label{lab14}
|A_{3}|\,\le\,\frac{1}{\tau_{0}}\sum\limits_{n>X}\frac{1}{\sqrt{n}}\biggl|\widehat{\Phi}\biggl(\frac{\ln{n}}{\tau_{0}}\biggr)\biggr|\,\le\,
\frac{1}{\tau_{0}}\sum\limits_{n>X}\frac{1}{\sqrt{n}}\exp{\biggl(-\,\frac{\ln{n}}{\tau_{0}}\,F\biggl(\frac{\ln{n}}{\tau_{0}}\biggr)\biggr)}\,=\\
=\,\frac{1}{\tau_{0}}\sum\limits_{m = 0}^{+\infty}\;\sum\limits_{Xe^{m\tau_{0}}<n\le Xe^{(m+1)\tau_{0}}}
\frac{1}{\sqrt{n}}\exp{\biggl(-\,\frac{\ln{n}}{\tau_{0}}\,F\biggl(\frac{\ln{n}}{\tau_{0}}\biggr)\biggr)}.
\end{multline}
Since $F$ is monotonic, we have
\[
F\biggl(\frac{\ln{n}}{\tau_{0}}\biggr)\,\ge\,F\biggl(\frac{1}{\tau_{0}}\,\ln{\bigl(Xe^{m\tau_{0}}\bigr)}\biggr)\,=\,F\biggl(\vf\biggl(\frac{\tau_{0}}{2}+1\biggr)+m\biggr)
\,\ge\,F\biggl(\vf\biggl(\frac{\tau_{0}}{2}+1\biggr)\biggr)\,=\,\frac{\tau_{0}}{2}+1
\]
for any $m\ge 0$ and for $Xe^{m\tau_{0}}<n\le Xe^{(m+1)\tau_{0}}$. Hence, the sum over $n$ in (\ref{lab14}) does not exceed
\begin{multline}\label{lab15}
\frac{1}{\tau_{0}}\sum\limits_{Xe^{m\tau_{0}}<n\le Xe^{(m+1)\tau_{0}}}
\frac{1}{\sqrt{n}}\exp{\biggl\{-\,\biggl(\vf\biggl(\frac{\tau_{0}}{2}+1\biggr)+m\biggr)\biggl(\frac{\tau_{0}}{2}+1\biggr)\biggr\}}\,\le\\
\le\, \frac{3}{\tau_{0}}\,\bigl(Xe^{(m+1)\tau_{0}}\bigr)^{0.5}\exp{\biggl\{-\,\biggl(\vf\biggl(\frac{\tau_{0}}{2}+1\biggr)+m\biggr)\biggl(\frac{\tau_{0}}{2}+1\biggr)\biggr\}}\,=\\
=\,\frac{3}{\tau_{0}}\,\exp{\biggl\{\frac{\tau_{0}}{2}\vf\biggl(\frac{\tau_{0}}{2}+1\biggr)\,+\,\frac{\tau_{0}}{2}\,(m+1)\,-\,
\biggl(\vf\biggl(\frac{\tau_{0}}{2}+1\biggr)+m\biggr)\biggl(\frac{\tau_{0}}{2}+1\biggr)\biggr\}}\,=\\
=\,\frac{3}{\tau_{0}}\,\exp{\biggl\{\frac{\tau_{0}}{2}\,-\,\vf\biggl(\frac{\tau_{0}}{2}+1\biggr)\biggr\}}e^{-\,m}\,<\,\frac{3e^{\,-m}}{\tau_{0}}.
\end{multline}
Finally we get
\begin{equation}\label{lab16}
|A_{3}|\,\le\,\frac{3}{\tau_{0}}\sum\limits_{m = 0}^{+\infty}e^{-m}\,<\,\frac{5}{\tau_{0}}.
\end{equation}
Since $\Phi(u)$ is non\,-negative, $|\widehat{\Phi}(\lambda)|\le
\widehat{\Phi}(0)$ for real $\lambda$. Hence, lemma 3 implies
\begin{multline}\label{lab17}
|A_{2}|\,\le\,\frac{1}{\tau_{0}}\sum\limits_{k\ge 2}\sum\limits_{p^{k}\le X}\frac{\widehat{\Phi}(0)}{kp^{0.5k}}\,\le
\,\frac{\widehat{\Phi}(0)}{\tau_{0}}\,\biggl(\frac{1}{2}\sum\limits_{p\le \sqrt{X}}\frac{1}{p}\,+\,\frac{1}{3}\sum\limits_{k\ge 3}\sum\limits_{p}p^{-0.5k}\biggr)\,<\\
<\,\frac{\widehat{\Phi}(0)}{2\tau_{0}}\,\biggl(\ln\ln{X}\,-\,\ln{2}\,+\,\mathfrak{m}\,+\,\frac{2}{3}\sum\limits_{p}\frac{1}{p(\sqrt{p}-1)}+\frac{4}{\ln^{2}{X}}\biggr)\,<\,
\frac{\widehat{\Phi}(0)}{2\tau_{0}}\,\biggl(\ln{\biggl(\tau_{0}\vf\biggl(\frac{\tau_{0}}{2}+1\biggr)\biggr)}\,+\,2\biggr)
\end{multline}
Summation of (\ref{lab13}), (\ref{lab16}), (\ref{lab17}) yields:
\begin{equation}\label{lab18}
|A_{2}|\,+\,|A_{3}|\,+\,|B|\,\le\,\frac{\widehat{\Phi}(0)}{\tau_{0}}\,\ln{\biggl(\tau_{0}\vf\biggl(\frac{\tau_{0}}{2}+1\biggr)\biggr)}.
\end{equation}
Now we set
\[
k\,=\,\biggl[\frac{e^{\alpha\tau_{0}}}{\alpha\kappa\tau_{0}}\biggr]\,\ge\,7,
\]
where $\kappa > 0$  will be chosen later. Denote
\[
a(p)\,=\,\widehat{\Phi}\biggl(\frac{\ln{p}}{\tau}\biggr),\quad V(t)\,=\,\sum\limits_{p\le X}\frac{a(p)}{\sqrt{p}}\,p^{it},\quad A_{0}(t)\,=\,0.5(V(t)\,+\,\overline{V}(t)),
\]
and define the integrals
\[
I(k)\,=\,\int_{T}^{T+H}A_{0}^{2k}(t)\,dt,\quad J(k)\,=\,\int_{T}^{T+H}A_{0}^{2k+1}(t)\,dt.
\]
Thus we find that
\[
I(k)\,=\,2^{-2k}\sum\limits_{\nu = 0}^{2k}\binom{2k}{\nu}j(\nu),\quad j(\nu)\,=\,\int_{T}^{T+H}V^{\nu}(t)\overline{V}^{\,\mu}(t)\,dt,\quad \mu = 2k-\nu.
\]
Setting for brevity $P = p_{1}\ldots p_{\nu}$, $Q = q_{1}\ldots
q_{\mu}$, in the case $\mu\ne\nu$ we get
\begin{multline*}
|j(\nu)|\,=\,\biggl|\sum\limits_{\substack{p_{1},\ldots, p_{\nu}\le
X \\ q_{1},\ldots, q_{\mu}\le X}}\frac{a(p_{1})\ldots
a(q_{\mu})}{\sqrt{p_{1}\ldots
q_{\mu}}}\int_{T}^{T+H}\biggl(\frac{p_{1}\ldots p_{\nu}}{q_{1}\ldots
q_{\mu}}\biggr)^{it}\,dt\biggr|\,\le\\
\le\,2\sum\limits_{P,Q}\frac{|a(p_{1})|\ldots
|a(q_{\mu})|}{\sqrt{PQ}}\biggl|\ln\frac{P}{Q}\biggr|^{-1}\,\le\,2\widehat{\Phi}^{2k}(0)\sum\limits_{P,Q}\frac{1}{\sqrt{PQ}}\biggl|\ln\frac{P}{Q}\biggr|^{-1}.
\end{multline*}
If $P<Q$ then
\[
\biggl|\ln\frac{P}{Q}\biggr|\,=\,\ln\frac{Q}{P}\,\ge\,\ln{\frac{P+1}{P}}\,\ge\,\frac{1}{2P};
\quad\text{otherwise,}\quad
\biggl|\ln\frac{P}{Q}\biggr|\,\ge\,\frac{1}{2Q}.
\]
Hence,
\begin{multline}\label{lab19}
|j(\nu)|\,\le\,2\widehat{\Phi}^{2k}(0)\biggl(\;\sum\limits_{P<Q}\frac{2P}{\sqrt{PQ}}\,+\,
\sum\limits_{P>Q}\frac{2Q}{\sqrt{PQ}}\biggr)\,\le\\
\le\,4\widehat{\Phi}^{2k}(0)\biggl(\;\sum\limits_{p_{1},\ldots,p_{\nu}\le
X}(p_{1}\ldots p_{\nu})^{0.5}\sum\limits_{q_{1},\ldots, q_{\mu}\le
X}(q_{1}\ldots
q_{\mu})^{-0.5}\,+\\
+\,\sum\limits_{p_{1},\ldots,p_{\nu}\le X}(p_{1}\ldots
p_{\nu})^{-0.5}\sum\limits_{q_{1},\ldots, q_{\mu}\le X}(q_{1}\ldots
q_{\mu})^{0.5}\biggr)\,=\\
=\,4\widehat{\Phi}^{2k}(0)\bigl(S^{\nu}C^{2k-\nu}\,+\,S^{2k-\nu}C^{\nu}\bigr),
\end{multline}
where
\[
S\,=\,\sum\limits_{p\le X}p^{0.5},\qquad C\,=\,\sum\limits_{q\le
X}q^{-0.5}.
\]
The same bound is true for the non-diagonal terms in the case
$\nu=k$. Summing (\ref{lab19}) over $0\le \nu\le 2k$ we get
\begin{equation}\label{lab20}
I(k)\,=\,2^{-2k}\binom{2k}{k}H\mathfrak{S}_{k}\,+\,8\theta
\widehat{\Phi}^{2k}(0)S^{2k},
\end{equation}
where
\[
\mathfrak{S}_{k}\,=\,\sum\limits_{p_{1}\ldots p_{k}=q_{1}\ldots
q_{k}}\frac{a^{2}(p_{1})\ldots a^{2}(p_{k})}{p_{1}\ldots
p_{k}},\quad |\theta|\le 1.
\]
Since
\[
S\,=\,\biggl(\frac{2}{3}\,+\,o(1)\biggr)\,\frac{X^{1.5}}{\ln{X}},
\]
then the last term in (\ref{lab20}) is less that
$(X^{1.5}(\ln{X})^{-0.5})^{2k}$ in modulus.

Estimating $\mathfrak{S}_{k}$ from below, we retain in
$\mathfrak{S}_{k}$ all the terms corresponding to the tuples
$(p_{1},\ldots,p_{k})$ without repetitions. Thus we get
\begin{equation}\label{lab21}
\mathfrak{S}_{k}\,\ge\,k!\sum\limits_{\substack{p_{1}\ldots
p_{k}=q_{1}\ldots q_{k} \\ p_{1},\ldots, p_{k}\;\text{are
distinct}}} \frac{a^{2}(p_{1})\ldots a^{2}(p_{k})}{p_{1}\ldots
p_{k}}.
\end{equation}
Since $\vf(v)$ is monotonic, we have $X\ge e^{\alpha\tau_{0}}$ for
sufficiently large $H$. Replacing the upper limit for
$p_{1},\ldots,p_{k}$ in (\ref{lab21}) by $e^{\alpha\tau_{0}}$ and
noting that
\[
\widehat{\Phi}\biggl(\frac{\ln{p}}{\tau_{0}}\biggr)\,\ge\,\widehat{\Phi}(\alpha)\,>\,0
\]
for $2\le p\le e^{\alpha\tau_{0}}$, we have
\begin{multline*}
\mathfrak{S}_{k}\,\ge\,k!\widehat{\Phi}^{2k}(\alpha)\sum\limits_{\substack{p_{1},\ldots,
p_{k}\le e^{\alpha\tau_{0}}\\ p_{1},\ldots, p_{k}\;\text{are
distinct}}}(p_{1}\ldots p_{k})^{-1}\,\ge\\
\ge\,k!\widehat{\Phi}^{2k}(\alpha)\sum\limits_{p_{1}\le
e^{\alpha\tau_{0}}}\frac{1}{p_{1}}\sum\limits_{\substack{p_{2}\le
e^{\alpha\tau_{0}}\\ p_{2}\ne p_{1}}}\frac{1}{p_{2}}\cdots
\sum\limits_{\substack{p_{k}\le e^{\alpha\tau_{0}}\\ p_{k}\ne
p_{1},\ldots,p_{k-1}}}\frac{1}{p_{k}}\,\ge\\
\ge\,k!\widehat{\Phi}^{2k}(\alpha)\sum\limits_{p_{1}\le
e^{\alpha\tau_{0}}}\frac{1}{p_{1}}\sum\limits_{\substack{p_{2}\le
e^{\alpha\tau_{0}}\\ p_{2}\ne \varpi_{1}}}\frac{1}{p_{2}}\cdots
\sum\limits_{\substack{p_{k}\le e^{\alpha\tau_{0}}\\ p_{k}\ne
\varpi_{1},\ldots,\varpi_{k-1}}}\frac{1}{p_{k}}\,\ge\\
\ge\,k!\widehat{\Phi}^{2k}(\alpha)\biggl(\;\sum\limits_{\varpi_{k-1}<p\le
e^{\alpha\tau_{0}}}\frac{1}{p}\biggr)^{\! k}.
\end{multline*}
Let us take $\kappa = \max{\bigl(61,4\kappa^{-1}\bigr)}$. Then, by
lemma 3 we get
\begin{multline*}
\varpi_{k-1}<k(\ln{k}+\ln\ln{k})\,\le\,\frac{e^{\alpha\tau_{0}}}{\alpha\kappa\tau_{0}}\,\bigl(\alpha\tau_{0}-\ln{(\alpha\kappa\tau_{0})}+\ln{(\alpha\tau_{0})}\bigr)
\,<\,\frac{e^{\alpha\tau_{0}}}{\kappa},\\
k\ln{k}\,>\,\biggl(\frac{e^{\alpha\tau_{0}}}{\alpha\kappa\tau_{0}}\,-\,1\biggr)\ln{\biggl(\frac{e^{\alpha\tau_{0}}}{\alpha\kappa\tau_{0}}\,-\,1\biggr)}\,>\,
\frac{e^{\alpha\tau_{0}}}{2\alpha\kappa\tau_{0}}\,\ln{\biggl(\frac{e^{\alpha\tau_{0}}}{2\alpha\kappa\tau_{0}}\biggr)}\,=\\
=\,\frac{e^{\alpha\tau_{0}}}{2\kappa}\,\biggl(1\,-\,\frac{\ln{(2\alpha\kappa\tau_{0})}}{\alpha\tau_{0}}\biggr),
\end{multline*}
and hence
\[
\ln{(k\ln{k})}\,>\,\alpha\tau_{0}\,\biggl(1\,-\,\frac{\ln{\kappa}+1}{\alpha\tau_{0}}\biggr),\quad
\frac{1}{\ln^{2\mathstrut}{(k\ln{k})}}\,<\,\frac{1}{(\alpha\tau_{0})^{2}}\,\biggl(1\,+\,\frac{3(\ln{\kappa}+1)}{\alpha\tau_{0}}\biggr).
\]
Using lemma 3 again, we obtain
\begin{multline*}
\sum\limits_{\varpi_{k-1}<p\le
e^{\alpha\tau_{0}}}\frac{1}{p}\,>\,\sum\limits_{k\ln{k}<p\le
e^{\alpha\tau_{0}}}\frac{1}{p}\,>\,\ln\ln{e^{\alpha\tau_{0}}}\,-\,\ln\ln{\frac{e^{\alpha\tau_{0}}}{\kappa}}\,
-\,\frac{1.5}{(\alpha\tau_{0})^{2}}\,\biggl(1\,+\,\frac{3(\ln{\kappa}+1)}{\alpha\tau_{0}}\biggr)\,=\\
=\,-\ln{\biggl(1\,-\,\frac{\ln{\kappa}}{\alpha\tau_{0}}\biggr)}\,-\,
\frac{1.5}{(\alpha\tau_{0})^{2}}\,\biggl(1\,+\,\frac{3(\ln{\kappa}+1)}{\alpha\tau_{0}}\biggr)\,=\\
=\,\frac{\ln{\kappa}}{\alpha\tau_{0}}\,+\,\frac{1}{2(\alpha\tau_{0})^{2}}\,\bigl((\ln{\kappa})^{2}\,-\,3\bigr)
\,+\,\frac{1}{3(\alpha\tau_{0})^{3}}\,\biggl((\ln{\kappa})^{3}\,-\frac{27}{2}(\ln\kappa\,+\,1)\biggr)\,>\,
\frac{\ln{\kappa}}{\alpha\tau_{0}}.
\end{multline*}
Passing to the estimation of $I(k)$ and noting that
$k=\D\biggl[\frac{1}{\alpha\kappa}\,\frac{\ln{H}}{\ln{X}}\biggr]$,
we find:
\begin{multline*}
\biggl(\frac{X^{1.5}}{\sqrt{\ln{X}}}\biggr)^{\!
2k}\,<\,X^{3k}\,\le\,\exp{\biggl(\frac{3\ln{H}}{\alpha\kappa}\biggr)}\,\le\,H^{0.75},\\
I(k)\,>\,2^{-2k}\binom{2k}{k}H
k!\widehat{\Phi}^{2k}(\alpha)\biggl(\frac{\ln{\kappa}}{\alpha\tau_{0}}\biggr)^{\!k}\,-\,\biggl(\frac{X^{1.5}}{\sqrt{\ln{X}}}\biggr)^{\!
2k}\,>\\
>\,\frac{(2k)!}{k!}\,H\,\biggl(\frac{\widehat{\Phi}(\alpha)}{4}\sqrt{\frac{\ln{\kappa}}{\alpha\tau_{0}}}\,\biggr)^{\!
2k}\,-\,H^{0.75}\,>\\
>\,\frac{e}{2}\biggl(\frac{4k}{e}\biggr)^{\!k}H\biggl(\frac{\widehat{\Phi}(\alpha)}{4}\sqrt{\frac{\ln{\kappa}}{\alpha\tau_{0}}}\,\biggr)^{\!
2k}\,-\,H^{0.75}\,>\,HM^{2k},
\end{multline*}
where
\[
M\,=\,\frac{\widehat{\Phi}(\alpha)}{2}\sqrt{\frac{k\ln{\kappa}}{e\alpha\tau_{0}}}\,>\,2.
\]

Repeating word-by-word the estimation of the non-diagonal terms of
$I(k)$, we get
\[
|J(k)|\,<\,\biggl(\frac{X^{1.5}}{\sqrt{\ln{X}}}\biggr)^{\!
2k+1}\,<\,X^{3k(1+1/(2k))}\,\le\,X^{4k}\,\le\,H^{4/(\alpha\kappa)}\,\le\,H,
\]
and hence $|J(k)|<0.5HM^{2k+1}$. By lemma 2, there exists $t_{0}$
such that $T\le t_{0}\le T+H$ and $A_{0}(t_{0})>0.5M$. Setting $t =
t_{0}$ in (\ref{lab11}) and taking into account (\ref{lab12}),
(\ref{lab18}) we find that
\begin{multline}\label{lab22}
I_{0}(t_{0})\,=\,A(t_{0})-B(t_{0})-I_{1}(t_{0})-I_{2}(t_{0})\,\ge\\
\ge\,\frac{M}{2\tau_{0}}\,-\,\frac{\widehat{\Phi}(0)}{\tau_{0}}\,\ln{\biggl(\tau_{0}\vf\biggl(\frac{\tau_{0}}{2}+1\biggr)\biggr)}
-\frac{4(\ln{t_{0}})}{\tau_{0}}\,\frac{e^{-\,G(\tau_{0}H)}}{G'(\tau_{0}H)}\,\ge\\
\ge\,\frac{\widehat{\Phi}(\alpha)}{4\tau_{0}}\biggl(\frac{e^{\alpha\tau_{0}}}{2\alpha\kappa\tau_{0}}\,\frac{\ln{\kappa}}{e\alpha\tau_{0}}\biggr)^{\!0.5}
\,-\,\frac{\widehat{\Phi}(0)}{\tau_{0}}\,\bigl(\ln{\tau_{0}}\,+\,e^{0.25\alpha\tau_{0}}\bigr)\,-\,\frac{4(\ln{t_{0}})e^{-\,\ln\ln{T}}}{\tau_{0}G'(\tau_{0}H)}
\,>\\
>\,\frac{\widehat{\Phi}(\alpha)}{10\alpha}\sqrt{\frac{\ln{\kappa}}{\kappa}}\,\frac{e^{0.5\alpha\tau_{0}}}{\tau_{0}^{2\mathstrut}}.
\end{multline}
The inequality (\ref{lab22}) and the definition of $I_{0}$ implies
that the maximum $M_{1}$ of the function
$\ln{|\zeta(0.5+i(t_{0}+u))|}$ on the segment $|u|\le H$ is strictly
positive. Hence,
\begin{equation}\label{lab23}
I_{0}(t_{0})\,\le\,M_{1}\int_{-H}^{H}\Phi(\tau_{0}H)\,du\,\le\,M_{1}\int_{-\infty}^{+\infty}\Phi(\tau_{0}H)\,du\,=\,\frac{\widehat{\Phi}(0)}{\tau_{0}}\,M_{1}.
\end{equation}
Comparing (\ref{lab22}) with (\ref{lab23}) and noting that the point
$t_{0}+u$ of maximum is contained in $\bigr[T-H,T+2H\bigl]$, we find
that
\begin{equation}\label{lab24}
\max_{T-H\le t\le
T+2H}\ln{\bigl|\zeta(0.5+it)\bigr|}\,\ge\,M_{1}\,>\,\frac{1}{10\alpha}\sqrt{\frac{\ln\kappa}{\kappa}}\,\frac{\widehat{\Phi}(\alpha)}{\widehat{\Phi}(0)}\,
\frac{e^{0.5\alpha\tau_{0}}}{\tau_{0}}.
\end{equation}
Thus, (\ref{lab07}) is proved.

Suppose now that $\vf(v)$ satisfies (\ref{lab08}). Then, setting
\[
X\,=\,e^{\Z \tau_{1}\vf\bigl(\frac{\Z \tau_{1}}{\Z
2\mathstrut}+1\bigr)},\quad
k\,=\,\biggl[\frac{e^{0.5\alpha\tau_{1}}}{\alpha\kappa\tau_{1}}\biggr]\,=\,\biggl[\frac{1}{\alpha\kappa}\,\frac{\ln{H}}{\ln{X}}\biggr],\quad
\kappa = \max{(0.5,4\alpha^{-1})}
\]
and repeating word-by-word the above arguments, we find that
\begin{multline*}
\varpi_{k-1}\,<\,e^{0.5\alpha\tau_{1}},\quad
\sum\limits_{\varpi_{k-1}<p\le
e^{\alpha\tau_{1}}}\,\ge\,\ln\ln{e^{\alpha\tau_{1}}}\,-\,\ln\ln{e^{0.5\alpha\tau_{1}}}\,-\,\frac{4.5}{(\alpha\tau_{1})^{2\mathstrut}}\,>\,\frac{2}{3},\\
\mathfrak{S}_{k}\,>\,2^{-2k}\,\frac{(2k)!}{k!}\,H\widehat{\Phi}^{2k}(\alpha)\,\biggl(\frac{2}{3}\biggr)^{k}\,-\,H^{0.75}\,>\,HM^{2k},
\end{multline*}
where $M = \widehat{\Phi}(\alpha)\sqrt{\frac{\D 2k}{\D 3e\mathstrut}}$,
and, similarly,
\[
|J(k)|\,<\,X^{15k/4}\,\le\,H^{15/16}\,<\,0.5HM^{2k+1}.
\]
By lemma 2, $A(t_{0})>0.5M$ for some $t_{0}$, $T\le t_{0}\le T+H$.
Since (\ref{lab08}) and (\ref{lab18}) imply the bound
\[
|A_{2}|\,+\,|A_{3}|\,+\,|B|\,\le\,\frac{\widehat{\Phi}(0)}{\tau_{1}}\,\ln{\biggl(\tau_{1}\vf\biggl(\frac{\tau_{1}}{2}+1\biggr)\biggr)}\,<\,0.5\widehat{\Phi}(0),
\]
we get
\[
\max_{T-H\le t\le
T+2H}\ln{|\zeta(0.5+it)|}\,>\,\frac{\widehat{\Phi}(\alpha)}{\sqrt{6e}}\,\sqrt{k}\,>\,\frac{\widehat{\Phi}(\alpha)}{6\sqrt{\alpha\kappa}}\,\frac{e^{0.25\alpha\tau_{1}}}{\sqrt{\tau_{1}}}.
\]
Theorem is proved.

\vspace{0.2cm}

To prove the Corollary, we use (\ref{lab24}) with $\alpha = \varrho
- \vep$, $\vep = 2\tau_{0}^{-1}$. For sufficiently large $H$ we have
$\widehat{\Phi}(\alpha)\,=\,-\widehat{\Phi}'(\varrho-\theta\vep)\vep\,\ge\,0.5\vep|\widehat{\Phi}'(\varrho)|\,=\,|\widehat{\Phi}'(\varrho)|\tau_{0}^{-1}$,
\begin{multline*}
\max_{T-H\le t\le
T+2H}\ln{|\zeta(0.5+it)|}\,>\,\frac{1}{10\varrho}\sqrt{\frac{\ln\kappa}{\kappa}}\,\frac{2}{\tau_{0}}\,\frac{|\widehat{\Phi}'(\varrho)|}{\widehat{\Phi}(0)}
\,\frac{1}{\tau_{0}}\,e^{0.5\tau_{0}(\varrho-2/\tau_{0})}\,>\\
>\,\frac{1}{5e\varrho}\sqrt{\frac{\ln{\kappa_{1}}}{\kappa_{1}}}\,\frac{|\widehat{\Phi}'(\varrho)|}{\widehat{\Phi}(0)}\,\frac{e^{0.5\tau_{0}\varrho}}{\tau_{0}^{2}},\quad
\end{multline*}
where $\kappa_{1}=\max{(62,5\varrho^{-1})}$. The second assertion of
the Corollary can be proved similarly.

\pagebreak

\textbf{\S 3. The rate of decreasing of some Fourier transforms}

\vspace{0.2cm}

In order to apply Theorem A for given function $\Phi(u)$, we need an
estimate of type (\ref{lab04}) for the rate of decreasing of
$\widehat{\Phi}(\lambda)$ when $\lambda\to\pm \infty$. In what follows,
we obtain some bounds of such type.

\vspace{0.2cm}

\textsc{Lemma 4.} \emph{Suppose $m\ge 1$ is any fixed integer,
$\Phi(u) = \exp{\biggl(-\,\frac{\D u^{2m}}{\D
2m\mathstrut}\biggr)}$. Then the inequality
\[
\bigl|\widehat{\Phi}(\lambda)\bigr|\,<\,\frac{5}{\sqrt{m}}\,|\lambda|^{-\beta}\,\exp{\biggl(-\,\frac{|\lambda|^{\alpha}}{\alpha}\,\sin\pi\kappa\biggr)}
\]
holds for any real $\lambda$, $|\lambda|\ge \lambda_{0}$, with}
\[
\alpha\,=\,\frac{2m}{2m-1},\quad \beta\,=\,\frac{m-1}{2m-1},\quad
\kappa\,=\,\frac{1}{2(2m-1)}.
\]

\vspace{0.2cm}

\textsc{Proof.} The case $m = 1$ is obvious. If $m\ge 2$, this
assertion follows from the asymptotic formula for
$\widehat{\Phi}(\lambda)$ from \cite[Ch. IV, \S 7]{Fedoryuk_1977}.

\vspace{0.2cm}

\textsc{Lemma 5.} \emph{Let $p,q$ be integers, $1\le p<q$,
$(p,q)=1$, $r = p/q$, $\vep = e^{\pi i/q}$, and let}
\[
G_{r}(z)\,=\,\sum\limits_{k =
0}^{q-1}\ch{\bigl(\vep^{k}z^{\,p/q}\bigr)},\quad
\Phi_{r}(z)\,=\,\exp{(-G_{r}(z))}.
\]
\emph{Then the estimate
$|\widehat{\Phi}_{r}(\lambda)|<\exp{\bigl(-|\lambda|\,F_{r}(|\lambda|)\bigr)}$
holds for any real $\lambda$, $|\lambda|>\lambda_{0}$, with}
\[
F_{r}(u)\,=\,\frac{3}{5}\biggl(\ln\frac{\lambda}{q}\biggr)^{\frac{\Z
q}{\Z p}-1}.
\]
\vspace{0.2cm}

\textsc{Proof.} Suppose that $\lambda>\lambda_{0}>0$ (the case of
negative $\lambda$ is treated in the same way). Since the function
\[
G_{r}(z)\,=\,q\sum\limits_{n=0}^{+\infty}\frac{z^{2np}}{(2nq)!}
\]
is entire function of order $r$, then, for any $y>0$, we have
\[
\widehat{\Phi}_{r}(\lambda)\,=\,\int_{-\infty}^{+\infty}e^{-G_{r}(z)-i\lambda
z}\,dz\,=\,\int_{-\infty}^{+\infty}e^{-G_{r}(x-iy)-i\lambda(x-iy)}dx\,=\,e^{-\lambda
y}\int_{-\infty}^{+\infty}e^{-G_{r}(x-iy)-i\lambda x}dx
\]
and hence
\[
|\widehat{\Phi}_{r}(\lambda)|\,\le\,e^{-\lambda
y}\int_{-\infty}^{+\infty}e^{-\RRe G_{r}(x-iy)}\,dx.
\]
In what follows, we suppose $y>y_{0}(p,q)$ to be sufficiently large
and set
\[
x-iy\,=\,\rho e^{-i\vf},\quad\text{where}\quad
\rho\,=\,\sqrt{x^{2}+y^{2}},\quad \vf\,=\,\arctg{\frac{y}{x}}.
\]
Then
\begin{multline}\label{lab25}
\RRe G_{r}(x-iy)\,=\,\RRe G_{r}(\rho
e^{-i\vf})\,=\,\RRe\sum\limits_{k=0}^{q-1}\ch{\bigl(\rho^{\frac{\Z
p}{\Z q}}e^{\frac{\Z i}{\Z q}\,(\pi k-p\vf)}\bigr)}\,=\\
=\,\sum\limits_{k=0}^{q-1}\ch{\biggl(\rho^{\frac{\Z p}{\Z
q}}\cos\frac{\pi k-p\vf}{q}\biggr)}\cos{\biggl(\rho^{\frac{\Z p}{\Z
q}}\sin\frac{\pi k-p\vf}{q}\biggr)}.
\end{multline}
Let
\[
x_{0}\,=\,\biggl(\frac{3\sqrt{2}}{\pi}\,y\biggr)^{\frac{\Z q}{\Z
q-p}}.
\]
Then, for $x\ge x_{0}$, we have
\begin{multline*}
x^{1\,-\,\frac{\Z p}{\Z q}}\,\ge\,\frac{3\sqrt{2}}{\pi}\,y,\qquad
0\le \vf\le \frac{y}{x}\,\le\,\frac{\pi}{3\sqrt{2}}\,x^{-\,\frac{\Z
p}{\Z q}},\\
\rho^{\frac{\Z p}{\Z
q}}\sin{\frac{p\vf}{q}}\,\le\,(x^{2}+y^{2})^{\frac{\Z p}{\Z
2q\mathstrut}}\,\frac{py}{qx}\,\le\,(2x^{2})^{\frac{\Z p}{\Z
2q\mathstrut}}\,\frac{\pi}{3\sqrt{2}}\,x^{-\,\frac{\Z p}{\Z
q}}\,<\,\frac{\pi}{3},
\end{multline*}
and hence
\begin{multline*}
\cos{\biggl(\rho^{\frac{\Z p}{\Z
q}}\sin{\frac{p\vf}{q}}\biggr)}\,>\,\cos\frac{\pi}{3}\,=\,\frac{1}{2},\\
\rho^{\frac{\Z p}{\Z q}}\cos{\frac{p\vf}{q}}\,\ge\,\rho^{\frac{\Z
p}{\Z q}}\,\biggl(1\,-\,\frac{1}{2}\biggl(\frac{p\vf}{q}\biggr)^{\!
2}\biggr)\,>\,x^{\frac{\Z p}{\Z
q}}\biggl(1\,-\,\frac{\vf^{2}}{2}\biggr)\,>\\
>\,x^{\frac{\Z p}{\Z
q}}\biggl(1\,-\,\frac{1}{2}\,x^{-\,\frac{\Z p}{\Z
2q\mathstrut}}\biggr)\,>\,x^{\frac{\Z p}{\Z q}}\,-\,x^{-\,\frac{\Z
p}{\Z q}}.
\end{multline*}
Denote by $A_{r}$ the term with $k=0$ in (\ref{lab25}). Then
\[
A_{r}\,=\,\ch{\biggl(\rho^{\frac{\Z p}{\Z q}}\cos{\frac{\vf
p}{q}}\biggr)}\cos{\biggl(\rho^{\frac{\Z p}{\Z q}}\sin{\frac{\vf
p}{q}}\biggr)}\,>\,\frac{1}{2}\,\ch{\bigl(x^{\frac{\Z p}{\Z
q}}\,-\,x^{-\,\frac{\Z p}{\Z
q}}\bigr)}\,>\,\frac{1}{5}\,\exp{\bigl(x^{\frac{\Z p}{\Z q}}\bigr)}.
\]
Suppose now that $1\le k\le q-1$. Then
\begin{multline*}
\biggl|\cos{\frac{\pi k-p\vf}{q}}\biggr|\,=\,\biggl|\cos{\frac{\pi k}{q}}\cos{\frac{p\vf}{q}}\,+\,\sin{\frac{\pi k}{q}}\sin{\frac{p\vf}{q}}\biggr|\,\le\\
\le\,\biggl|\cos{\frac{\pi k}{q}}\biggr|\,+\,\vf\,<\,\biggl|\cos{\frac{\pi k}{q}}\biggr|\,+\,\frac{\pi}{3\sqrt{2}}\,x^{-\,\frac{\Z p}{\Z q}}.
\end{multline*}
Since
\[
\rho^{\frac{\Z p}{\Z q}}\,=\,x^{\frac{\Z p}{\Z q}}\biggl(1\,+\,\frac{y^{2}}{x^{2\mathstrut}}\biggr)^{\frac{\Z p}{\Z 2q\mathstrut}}\,<\,
x^{\frac{\Z p}{\Z q}}\bigl(1\,+\,x^{-\,\frac{\Z p}{\Z 2q\mathstrut}}\bigr)^{\frac{\Z p}{\Z 2q\mathstrut}}\,<\,x^{\frac{\Z p}{\Z q}}\biggl(1\,+\,\frac{p}{2q}\,x^{-\,\frac{\Z p}{\Z 2q\mathstrut}}\biggr)\,<\,x^{\frac{\Z p}{\Z q}}\,+\,0.5\,x^{-\,\frac{\Z p}{\Z q}},
\]
we get
\[
\biggl|\rho^{\frac{\Z p}{\Z q}}\cos{\frac{\pi k-p\vf}{q}}\biggr|\,<\,\bigl(x^{\frac{\Z p}{\Z q}}\,+\,0.5\,x^{-\,\frac{\Z p}{\Z q}}\bigr)
\biggl(\biggl|\cos{\frac{\pi k}{q}}\biggr|\,+\,\frac{\pi}{3\sqrt{2}}\,x^{-\,\frac{\Z p}{\Z q}}\biggr)\,<\,
x^{\frac{\Z p}{\Z q}}\biggl|\cos{\frac{\pi k}{q}}\biggr|\,+\,\frac{3}{4}.
\]
Thus,
\[
\ch{\biggl(\rho^{\frac{\Z p}{\Z
q}}\cos{\frac{p\vf}{q}}\biggr)}\,\le\,\ch{\biggl(x^{\frac{\Z p}{\Z
q}}\biggl|\cos{\frac{\pi k}{q}}\biggr|\,+\,\frac{3}{4}\biggr)}
\,<\,1.06\exp{\biggl(x^{\frac{\Z p}{\Z q}}\biggl|\cos{\frac{\pi
k}{q}}\biggr|\biggr)}.
\]
Denote by $B_{r}$ the sum in (\ref{lab25}) of the terms with $k>0$.
Then
\[
|B_{r}|\,<\,1.06\sum\limits_{k=1}^{q-1}\exp{\biggl(x^{\frac{\Z p}{\Z
q}}\biggl|\cos{\frac{\pi
k}{q}}\biggr|\biggr)}\,<\,2.2\exp{\biggl(x^{\frac{\Z p}{\Z
q}}\cos{\frac{\pi}{q}}\biggr)}.
\]
Hence,
\[
\RRe
G_{r}(x-iy)\,\ge\,A_{r}-|B_{r}|\,>\,\frac{1}{5}\,\exp{\bigl(x^{\frac{\Z
p}{\Z q}}\bigr)}\,-\,\frac{11}{5}\,\exp{\biggl(x^{\frac{\Z p}{\Z
q}}\cos{\frac{\pi}{q}}\biggr)}\,>\,\frac{1}{6}\exp{\bigl(x^{\frac{\Z
p}{\Z q}}\bigr)}
\]
for any $x\ge x_{0}$. The similar bound (with $|x|$ instead of $x$)
holds for $x\le -x_{0}$. If $|x|\le x_{0}$ then
\begin{multline*}
\rho\,\le\,(x_{0}^{2}+y^{2})^{0.5}\,=\,\biggl\{\biggl(\frac{3\sqrt{2}}{\pi}\,y\biggr)^{\frac{\Z
2p}{\Z q-p\mathstrut}}\,+\,y^{2}\biggr\}^{0.5}\,=\\
=\,y^{\frac{\Z q}{\Z
q-p\mathstrut}}\biggl\{\biggl(\frac{3\sqrt{2}}{\pi}\biggr)^{\frac{\Z
2q}{\Z q-p\mathstrut}}\,+\,y^{-\,\frac{\Z 2p}{\Z
q-p\mathstrut}}\biggr\}^{0.5}\,<\,\bigr(y\sqrt{2}\bigl)^{\frac{\Z
q}{\Z q-p\mathstrut}},
\end{multline*}
so we have
\[
\bigl|\RRe G_{r}(x-iy)\bigr|\,\le\,q\ch{\bigl(\rho^{\frac{\Z p}{\Z
q}}\bigr)}\,<\,q\ch{\bigl((y\sqrt{2})^{\frac{\Z p}{\Z q-p}}\bigr)}.
\]
Passing to the estimate of $\widehat{\Phi}(\lambda)$, we obtain
\begin{multline*}
\bigl|\widehat{\Phi}(\lambda)\bigr|\,\le\,e^{-\lambda
y}\biggl(\;\int_{-x_{0}}^{x_{0}}\exp{\biggl\{q\ch{\bigl((y\sqrt{2})^{\frac{\Z
p}{\Z q-p}}\bigr)}\biggr\}}dx \,+\\
+\,\biggl(\;\int_{-\infty}^{-x_{0}}\,+\,\int_{x_{0}}^{+\infty}\biggr)\exp{\biggl(-\frac{1}{6}\exp{\bigl(x_{0}^{\frac{\Z
p}{\Z q}}\bigr)}\biggr)}\biggr)dx\,<\\
<\,e^{-\lambda
y}\biggl(2x_{0}\exp{\biggl\{q\ch{\bigl((y\sqrt{2})^{\frac{\Z p}{\Z
q-p}}\bigr)}\biggr\}}\,+\,\exp{\biggl\{-\frac{1}{6}\exp{\bigl(x_{0}^{\frac{\Z
p}{\Z q}}\bigr)}\biggr\}}\biggr)\biggr).
\end{multline*}
Since $\sqrt{2}<1.5$, we have
\[
\bigl|\widehat{\Phi}(\lambda)\bigr|\,<\,\exp{\biggl(-\,\lambda
y\,+\,0.5qe^{(1.5y)^{\frac{\Z p}{\Z q-p}}}\biggr)}.
\]
Setting $y = \frac{\D 2}{\D
3\mathstrut}\bigl(\ln(\lambda/q)\bigr)^{\frac{\Z q}{\Z p}-1}$, we
finally get
\[
\bigl|\widehat{\Phi}(\lambda)\bigr|\,<\,\exp{\biggl(-\,\lambda y + 0.5q
e^{\ln\frac{\Z \lambda}{\Z
q}}\biggr)}\,=\,\exp{\bigl(-\lambda(y-0.5)\bigr)}\,<\,\exp{\biggl(-\,\frac{3}{5}\lambda\,\biggl(\ln\frac{\lambda}{q}\biggr)^{\frac{\Z
q}{\Z p}-1}\biggr)}.
\]
Lemma is proved.

\vspace{0.2cm}

\textsc{Lemma 6.} \emph{Let $G(u) = G_{1/2}(u)$, $\Phi(u)
=\Phi_{1/2}(u)$} ( \emph{in notations of lemma 5}). \emph{Then the
inequality}
\[
\bigl|\widehat{\Phi}(\lambda)\bigr|\,\le\,\exp{\biggl(-\,\frac{\pi}{1+\delta}\,|\lambda|\ln{|\lambda|}\biggr)}
\]
\emph{holds for any fixed $\delta$, $0<\delta<\delta_{0}<0.5$, and
for any real} $\lambda$, $|\lambda|> \lambda_{0}(\delta)$.

\vspace{0.2cm}

\textsc{Proof.} Applying the same arguments as above, we get
\[
\bigl|\widehat{\Phi}(\lambda)\bigr|\,\le\,e^{-\lambda
y}\int_{-\infty}^{+\infty}e^{-\RRe G(x-iy)}dx.
\]
Let $\delta_{1}, \delta_{2}, \ldots$ denote some positive constants
depending on $\delta$ and such that $\delta_{j}\to 0$ when
$\delta\to 0$. Taking $\rho = \sqrt{x^{2}+y^{2}}$, $\varphi =
\arctg{(y/x)}$, $x_{0} = (1+\delta_{1})(y/\pi)^{2}$, we get for
$x\ge x_{0}$:
\begin{multline*}
0<\sqrt{\rho}\sin\frac{\vf}{2}\,\le\,\sqrt{x}\biggl(1\,+\,\biggl(\frac{y}{x}\biggr)^{2}\biggr)^{0.25}\,\frac{y}{2x}\,\le\,\frac{y}{2\sqrt{x}}\,
(1+\delta_{2})\,<\,\frac{\pi}{2}\,(1-\delta),\\
0<\varphi < \frac{y}{x}\,<\,\frac{\pi}{\sqrt{x}}.
\end{multline*}
Hence,
\begin{multline*}
\sqrt{\rho}\cos\frac{\vf}{2}\,>\,\sqrt{x}\biggl(1\,-\,\frac{\vf^{2}}{8}\biggr)\,>\,\sqrt{x}\biggl(1\,-\,\frac{\pi^{2}}{8x}\biggr)\,>\,\sqrt{x}
\,-\,\frac{5}{4\sqrt{x}},\\
\ch{\biggl(\sqrt{\rho}\cos\frac{\vf}{2}\biggr)}\cos{\biggl(\sqrt{\rho}\sin\frac{\vf}{2}\biggr)}\,>\,
\frac{1}{2}e^{\sqrt{x} \,-\,\frac{\Z 5}{\Z
4\sqrt{x}\mathstrut}}\cos{\frac{\pi}{2}(1-\delta)}\,=\\
=\,\frac{1}{2}e^{\sqrt{x} \,-\,\frac{\Z 5}{\Z
4\sqrt{x}\mathstrut}}\sin\frac{\pi\delta}{2},\\
\biggl|\ch{\biggl(\sqrt{\rho}\sin\frac{\vf}{2}\biggr)}\cos{\biggl(\sqrt{\rho}\cos\frac{\vf}{2}\biggr)}\biggr|\,
\le\,\ch{\biggl(\sqrt{\rho}\sin\frac{\vf}{2}\biggr)}\,<\,\ch\frac{\pi}{2}.
\end{multline*}
Thus, we have
\begin{multline*}
\RRe G(x-iy)
=\ch{\biggl(\sqrt{\rho}\cos{\frac{\vf}{2}}\biggr)}\cos{\biggl(\sqrt{\rho}\sin{\frac{\vf}{2}}\biggr)}\,+\,
\ch{\biggl(\sqrt{\rho}\sin{\frac{\vf}{2}}\biggr)}\cos{\biggl(\sqrt{\rho}\cos{\frac{\vf}{2}}\biggr)}\,>\\
\,>\,\ch{\biggl(\sqrt{\rho}\cos{\frac{\vf}{2}}\biggr)}\cos{\biggl(\sqrt{\rho}\sin{\frac{\vf}{2}}\biggr)}\,-\,\ch{\biggl(\sqrt{\rho}\sin{\frac{\vf}{2}}\biggr)}\,>\\
>\,\frac{1}{2}\biggl(\sin\frac{\pi\delta}{2}\biggr)e^{\sqrt{x}\,-\,\frac{\Z
5}{\Z
4\sqrt{x}}}\,-\,\ch\frac{\pi}{2}\,>\,\frac{3\delta}{4}\,e^{\sqrt{x}}.
\end{multline*}
for $x$ under considering. The same bound is true for $x\le -x_{0}$
with $|x|$ instead of $x$. In the case $|x|\le x_{0}$, we have
\[
\sqrt{\rho}\,\le\,(x_{0}^{2}\,+\,y^{2})^{0.25}\,\le\,\biggl((1+\delta_{1})^{2}\biggl(\frac{y}{\pi}\biggr)^{4}\,+\,y^{2}\biggr)^{0.25}\,<\,
(1+\delta_{3})\,\frac{y}{\pi},
\]
and hence
\[
|\RRe
G(x-iy)|\,\le\,2\ch{\bigl(\sqrt{\rho}\bigr)}\,<\,\exp{\biggl((1+\delta_{4})\,\frac{y}{\pi}\biggr)}.
\]
Thus we obtain
\begin{multline*}
\bigl|\widehat{\Phi}(\lambda)\bigr|\,\le\,e^{-\lambda
y}\biggl(\;\int_{-x_{0}}^{x_{0}}\exp{\biggl\{\exp{\biggl((1+\delta_{4})\,\frac{y}{\pi}\biggr)}\biggr\}}\,dx\,+\\
+\,2\int_{x_{0}}^{+\infty}\exp{\biggl\{-\,\frac{3\delta}{4}\,e^{\sqrt{x}}\biggr\}}dx\biggr)\,
<\,3x_{0}\exp{\left\{-\lambda y\,+\,e^{(1+\delta_{4})\,\frac{\Z
y}{\Z \pi}}\right\}}.
\end{multline*}
Now let $y = \frac{\D \pi\ln{\lambda}}{\D 1+\delta_{4}\mathstrut}$.
Then
\[
\bigl|\widehat{\Phi}(\lambda)\bigr|\,<\,\frac{3(1+\delta_{1})}{(1+\delta_{4})^{2}}\,
(\ln{\lambda})^{2}\,\exp{\biggl\{-\,\frac{\pi\lambda\ln\lambda}{1+\delta_{4}}\,+\,\lambda\biggr\}}.
\]
If $\delta$ is sufficiently small, then
\[
\bigl|\widehat{\Phi}(\lambda)\bigr|\,<\,\exp{\biggl\{-\,\frac{\pi\lambda\ln\lambda}{1+\delta}\biggr\}}.
\]
The case of negative $\lambda$ can be treated in the same way. Lemma
is proved.

\vspace{0.2cm}

\textbf{\S 4. Basic assertions}

\vspace{0.2cm}

Here we prove Theorems 1\,-4. In what follows, we use the notations
of \S 2, 3 without any comments.

\vspace{0.2cm}

\textsc{Proof of Theorem 1.} Let $m\ge 2$, $\Phi(u) =
\exp{\biggl(-\,\frac{\D u^{2m}}{\D 2m\mathstrut}\biggr)}$. By lemma
4, the estimate (\ref{lab04}) holds for
\[
F(\lambda)\,=\,c_{0}\lambda^{\frac{\Z 1}{\Z 2m-1\mathstrut}},\quad
c_{0}\,=\,\frac{\sin{\pi\kappa}}{1+2\kappa},\quad
\kappa\,=\,\frac{1}{2(2m-1)}
\]
and for sufficiently large $|\lambda|$. Obviously, we have $\vf(v) =
(cv)^{2m-1}$, $c = c_{0}^{-1}$. Hence, the equation (\ref{lab09})
takes the form
\[
\frac{\alpha\tau_{1}}{2}\,+\,(2m-1)\biggl(\ln{\biggl(\frac{\tau_{1}}{2}+1\biggr)}\,+\,\ln{c}\biggr)\,=\,\ln\ln{H}.
\]
For fixed $m, \alpha$ and $H\to +\infty$, we have
\[
\frac{\alpha\tau_{1}}{2}\,=\,\ln\ln{H}\,-\,(2m-1)\ln\ln\ln{H}\,+\,(2m-1)
\ln{\alpha c_{0}}\,+\,O\biggl(\frac{\ln\ln\ln{H}}{\ln\ln{H}}\biggr),
\]
so hence
\[
\frac{e^{0.25\alpha\tau_{1}}}{\sqrt{\tau_{1}}}\,>\,(\alpha
c_{0})^{m}\,\frac{\sqrt{\ln{H}}}{(\ln\ln{H})^{m}}\,\ge\,\biggl(\frac{2\alpha\kappa}{2\kappa+1}\biggr)^{m}\,\frac{\sqrt{\ln{H}}}{(\ln\ln{H})^{m}}\,=\,
\frac{\alpha^{m}\sqrt{\ln{H}}}{(2m\ln\ln{H})^{m}}.
\]
Since $g(v) = (2mv)^{1/(2m)}$, then for any $H\ge
(1/3)(2m\ln\ln{T})^{1/(2m)}$ and some $\alpha>0$ we obtain from
(\ref{lab10}) that
\[
F(T;H)\,>\,\exp{\biggl(\frac{1}{6\sqrt{\alpha\kappa}}\,\frac{\widehat{\Phi}(\alpha)}{\widehat{\Phi}(0)}\,\frac{\alpha^{m}\sqrt{\ln{H}}}{(2m\ln\ln{H})^{m}}\biggr)},
\quad\text{where}\quad \kappa\,=\,\max{(0.5,4\alpha^{-1})}.
\]
One can check (see \cite{Korolev_2013}) that $\widehat{\Phi}(u)$ is
positive and monotonically decreasing for $0\le u\le 1$ and
\[
\widehat{\Phi}(0)\,=\,2(2m)^{\frac{\Z 1}{\Z
2m\mathstrut}}\Gamma\biggl(1\,+\,\frac{1}{2m}\biggr),\quad
\widehat{\Phi}(1)\,>\,\frac{5}{4}\,\exp{\biggl(-\,\frac{1}{2m}\biggl(\frac{\pi}{4}\biggr)^{\!
2m}\biggr)}.
\]
Finally we get
\begin{multline*}
F(T;H)\,>\,\exp{\biggl\{\frac{5}{96}\,\exp{\biggl(-\,\frac{1}{2m}\biggl(\frac{\pi}{4}\biggr)^{\!2m}\biggr)}\,(2m)^{-\,\frac{\Z
1}{\Z
2m\mathstrut}}\Gamma^{\,-\,1}\biggl(1\,+\,\frac{1}{2m}\biggr)\,
\frac{\sqrt{\ln{H}}}{(2m\ln\ln{H})^{m}}\biggr\}}\,>\\
>\,\exp{\biggl(\frac{0.05\sqrt{\ln{H}}}{(2m\ln\ln{H})^{m}}\biggr)}.
\end{multline*}
Theorem 1 is proved.

\vspace{0.2cm}

\textsc{Proof of Theorem 2.} Let $r = p/q<0.5$, $\Phi(u) =
\Phi_{r}(u)$. By lemma 5, one can take
\[
F(\lambda)\,=\,\frac{3}{5}\biggl(\ln\frac{\lambda}{q}\biggr)^{\frac{\Z
q}{\Z p}-1},\qquad
\vf(v)\,=\,q\exp{\biggl(\biggl(\frac{5v}{3}\biggr)^{\frac{\Z p}{\Z
q-p}}\biggr)}
\]
to satisfy (\ref{lab04}). Thus, (\ref{lab06}) takes the form
\begin{equation}\label{lab26}
\alpha\tau_{0}\,+\,\biggl(\frac{5}{3}\biggl(\frac{\tau_{0}}{2}+1\biggr)\biggr)^{\frac{\Z
p}{\Z q-p}}\,+\,\ln{q}\,=\,\ln\ln{H}.
\end{equation}
Since $0<c=\frac{\D p}{\D q-p\mathstrut}<1$, the solution $\tau_{0}$
satisfies the relation
\[
\alpha\tau_{0}\,=\,\ln\ln{H}\,-\,\biggl(\frac{5}{6\alpha}\ln\ln{H}\biggr)^{c}\,-\,\ln{q}\,+\,O\bigl((\ln\ln{H})^{2c-1}\bigr).
\]
One can check that
\[
\frac{e^{0.5\alpha\tau_{0}}}{\tau_{0}}\,>\,\alpha\sqrt{\frac{\ln
H}{q}}\exp{\biggl\{-\frac{1}{2}\,\biggl(\frac{5}{6\alpha}\,\ln\ln{H}\biggr)^{c}+O\bigl((\ln\ln{H})^{2c-1}\bigr)\biggr\}}(\ln\ln{H})^{-1}.
\]
Since $g(v) = (\ln{(3v)})^{q/p}$, we have for $H\,\ge\,
(1/3)(\ln\ln\ln{T})^{q/p}$:
\[
F(T;H)\,>\,\exp{\biggl(\sqrt{\ln{H}}\,e^{-(c_{0}\ln\ln{H})^{\frac{\Z
p}{\Z q-p}}}\biggr)},\qquad c_{0} = \frac{1}{\alpha}.
\]
In particular, for $q = 2m+1$, $p=m$ and $H\,\ge\,
(\ln\ln\ln{T})^{2+1/m}$ we have
\[
F(T;H)\,>\,\exp{\biggl(\sqrt{\ln{H}}\,e^{-c_{1}(\ln\ln{H})^{1-\frac{\Z
1}{\Z m+1\mathstrut}}}\biggr)},\quad c_{1} = c_{1}(m).
\]
Given $\vep$, we define $m$ by the conditions $\frac{\D 1}{\D m}\le
\vep < \frac{\D 1}{\D m-1\mathstrut}$. Then for any $H$,
\[
H\,>\,\frac{1}{3}\,(\ln\ln\ln{T})^{2+\vep}\,\ge\,\frac{1}{3}(\ln\ln\ln{T})^{2+\frac{\Z
1}{\Z m}}
\]
we obtain:
\[
F(T;H)\,>\,\exp{\biggl(\sqrt{\ln{H}}\,e^{-c_{1}(\ln\ln{H})^{1-\frac{\Z
1}{\Z
m+1\mathstrut}}}\biggr)}\,>\,\exp{\biggl(\sqrt{\ln{H}}\,e^{-c_{2}(\ln\ln{H})^{\Z
1-0.5\vep}}\biggr)}
\]
for some $c_{2} = c_{2}(\vep)>0$. Theorem 2 is proved.

\vspace{0.2cm}

\textsc{Proof of Theorem 3.} Let $\Phi(u) = \Phi_{1/2}(u) = e^{-(\ch{\sqrt{u}}+\cos{\sqrt{u}})}$. Since
$|\Phi(u)|\le e^{1-0.5e^{\sqrt{\Z |u|}}}$ for real $u$, one can take
\[
G(u)\,=\,\frac{1}{2}\,e^{\sqrt{u}}\,-\,1,\qquad g(v)\,=\,\ln^{2}(2v+2).
\]
Given $\delta>0$, lemma 5 implies that
\[
\bigl|\widehat{\Phi}(\lambda)\bigr|\,\le\,\exp{\biggl(-\;\frac{\pi|\lambda|\ln{|\lambda|}}{1+\delta}\biggr)}
\]
for any real $\lambda$, $|\lambda|>\lambda_{0}(\delta)$. Therefore,
\[
\vf(v)\,=\,\exp{\biggl((1+\delta)\,\frac{v}{\pi}\biggr)},
\]
and (\ref{lab06}) takes the form
\[
\alpha\tau_{0}\,+\,\frac{1+\delta}{\pi}\,\biggl(\frac{\tau_{0}}{2}\,+\,1\biggr)\,=\,\ln\ln{H}.
\]
Hence,
\[
\tau_{0}\,=\,\frac{\ln\ln{H}-\,\frac{\D 1+\delta}{\D \pi}}{\alpha\,+\,\frac{\D 1+\delta}{\D 2\pi\mathstrut}}.
\]
Let $\alpha = \varrho_{1}$ be the least positive root of the function $\widehat{\Phi}(\lambda)$. Then
\begin{multline*}
e^{0.5\varrho_{1}\tau_{0}}\,=\,\exp{\biggl\{\frac{0.5\varrho_{1}}{\varrho_{1}+(1+\delta)/(2\pi)}\,\ln\ln{H}\,-\,\frac{\varrho_{1}(1+\delta)}{2\pi\varrho_{1}+1+\delta}\biggr\}}\,>\\
>\,\exp{\biggl\{\frac{1-\delta}{2+(\pi\varrho_{1})^{-1\mathstrut}}\,\ln\ln{H}\biggr\}}\,=\,(\ln{H})^{\frac{\Z 1-\delta}{\Z 2+(\pi\varrho_{1})^{-1\mathstrut}}}.
\end{multline*}
Given $\vep>0$, we can choose $\delta$ to satisfy the inequalities
\[
F(T;H)\,>\,\exp{\biggl\{\frac{1}{5e\varrho_{1}}\sqrt{\frac{5\ln{2}}{8e}}\,\frac{|\widehat{\Phi}'(\varrho_{1})|}{\widehat{\Phi}(0)}\,\frac{e^{0.5\varrho_{1}\tau_{0}}}{\tau_{0}^{2}}\biggr\}}\,>\,
\exp{\bigl((\ln{H})^{\gamma-\vep}\bigr)}
\]
for any $H \ge (\ln\ln\ln{T})^{2}$ and $\gamma = \frac{\D 1}{\D 2 + (\pi\varrho_{1})^{-1\mathstrut}}$. The approximate calculations in ``Wolfram Mathematica 7.0''
show that $2.37689234<\varrho_{1}<2.37689235$. Hence, $\gamma = 0.46862145\ldots$. Theorem 3 is proved.

\vspace{0.2cm}

\textsc{Proof of Theorem 4.} Let $0.5<r = p/q<1$, $\Phi(u) =
\Phi_{r}(u)$. Similarly to the proof of Theorem 2, one can check that the equation (\ref{lab06}) has the form
(\ref{lab26}). Since $c = \frac{\D p}{\D q-p\mathstrut}>1$, its solution satisfies the relation
\[
\tau_{0}\,=\,\frac{6}{5}(\ln\ln{H})^{\frac{\Z 1}{\Z
c\mathstrut}}\,-\,\frac{36\alpha}{25 c}(\ln\ln{H})^{\frac{\Z 2}{\Z
c\mathstrut}-1}\,-2\,+\,O\bigl((\ln\ln{H})^{\eta}\bigr),
\]
where $\eta = \min{\biggl(\frac{\D 3}{\D c}-2,\frac{\D 1}{\D
c}-1\biggr)}$. Hence, for $H\,\ge\,(\ln\ln\ln{T})^{\frac{\Z q}{\Z
p}}$, we have
\begin{align*}
&
\frac{e^{0.5\alpha\tau_{0}}}{\tau_{0}}\,>\,\exp{\biggl\{\frac{3\alpha}{5}(\ln\ln{H})^{\frac{\Z
1}{\Z c}}\,-\,\frac{18\alpha^{2}}{25c}(\ln\ln{H})^{\frac{\Z 2}{\Z
c}-1}
\,-\,\alpha\,+\,O\bigl((\ln\ln{H})^{\eta}\bigr)\biggr\}}(\ln\ln{H})^{-\,\frac{\Z
1}{\Z c}}, \\
& F(T;H)\,>\,\exp{\bigl\{\exp{\bigl(0.5\alpha(\ln\ln{H})^{\frac{\Z
q}{\Z p}-1}\bigr)}\bigr\}}.
\end{align*}
Given $\vep$, we define $m$ by the inequalities $\frac{\D 1}{\D
m}\le \vep\!<\!\frac{\D 1}{\D m-1\mathstrut}$ and set $q = m+1, p =
m$. Taking $H\,\ge\,(\ln\ln\ln{T})^{1+\vep}\,\ge\,
(\ln\ln\ln{T})^{1+\frac{\Z 1}{\Z m\mathstrut}}$, we obtain:
\[
F(T;H)\,>\,\exp{\bigl\{\exp{\bigl(0.5\alpha(\ln\ln{H})^{\frac{\Z
1}{\Z m}}\bigr)}\bigr\}}\,>\,
\exp{\bigl\{\exp{\bigl(0.5(\ln\ln{H})^{0.5\vep}\bigr)}\bigr\}}.
\]
Theorem 4 is proved.

\vspace{0.5cm}

\renewcommand{\refname}{\normalsize{Bibliography}}


\begin{thebibliography}{99}

\bibitem{Balasubramanian_1986}

R.~Balasubramanian, ``On the frequency of Titchmarsh's phenomenon
for $\zeta(s)$. IV'', \emph{Hardy\,-Ramanujan J}., \textbf{9}(1986),
1\,-10.

\bibitem{Farmer_Gonek_Hughes_2007}

D.W.~Farmer, S.M.~Gonek, C.P.~Hughes, ``The maximum size of
$L$\,-functions'', \emph{J. Reine angew. Math.}, \textbf{609}(2007),
215\,-236.

\bibitem{Karatsuba_2001a}

A.A.~Karatsuba, On lower bounds for the Riemann zeta function,
\emph{Doklady Mathematics}, \textbf{63}:1 (2001), 9\,-10.

\bibitem{Karatsuba_2001b}

A.A.~Karatsuba, ``Lower bounds for the maximum modulus of $\zeta(s)$
in small domains of the critical Strip'', \emph{Math. Notes},
\textbf{70}:5 (2001), 724\,–726.

\bibitem{Karatsuba_2004}

A.A.~Karatsuba, ``Lower bounds for the maximum modulus of the
Riemann zeta function on short segments of the critical line''
\emph{Izv. Math.}, \textbf{68}:6 (2004), 1157\,–1163.

\bibitem{Karatsuba_2006}

A.A.~Karatsuba, ``Zero multiplicity and lower bound estimates of
$|\zeta(s)|$, \emph{Funct. Approx. Comment. Math.}, \textbf{35}
(2006), 195\,–207.

\bibitem{Garaev_2002}

M.Z.~Garaev, Concerning the Karatsuba conjectures, \emph{Taiwanese
J. Math.}, \textbf{6}:4 (2002), 573.

\bibitem{Feng_2004}

S.J.~Feng, ``On Karatsuba conjecture and the Lindelof hypothesis'',
\emph{Acta Arithmetica},\textbf{114}:3 (2004), 295.

\bibitem{Changa_2004}

M.E.~Changa, ``Lower Bounds for the Riemann Zeta Function on the
Critical Line'', \emph{Math. Notes}, \textbf{76}:5-6 (2004),
859\,-864.

\bibitem{Changa_2005}

M.E.~Changa, ``On a function-theoretic inequality'', \emph{Russian
Math. Surveys}, 2005, \textbf{60}:3 (2005), 564\,–565.

\bibitem{Korolev_2014}

M.A.~Korolev, ``On large values of the Riemann zeta-function on
short segments of the critical line'', Acta Arithm., 166 (2014), №
4, 349\,-390.

\bibitem{Selberg_1946}

A.~Selberg, ``Contributions to the theory of the Riemann
zeta-function'', \emph{Archiv Math. Naturvid.}, \textbf{48}:5(1946),
89-155 (see also: A.~Selberg, Collected papers. Vol. I. -- Berlin,
Springer\,-Verlag, 1989, 214\,-280).

\bibitem{Tsang_1984}

K.-M.~Tsang, The distribution of the values of the Riemann
zeta\,-function. A dissertation presented to the Faculty of
Princeton University in candidacy for the degree of Doctor of
Philosofy. Princeton, October 1984.

\bibitem{Tsang_1986}

K.-M.~Tsang, ``Some $\Omega$\,-theorems for the Riemann
zeta-function'', \emph{Acta Arith.}, \textbf{46}(1986), 369\,-395.

\bibitem{Tsang_1993}

K.-M.~Tsang, ``The large values of the Riemann zeta-function'',
\emph{Mathematika}, \textbf{40}(1993), 203\,-214.

\bibitem{Korolev_2005}

M.A.~Korolev, ``On large values of the function $S(t)$ on short
intervals'', \emph{Izv. Math.}, \textbf{69}:1 (2005), 113\,–122.

\bibitem{Boyarinov_2011a}

R.N.~Boyarinov, ``Omega theorems in the theory of the Riemann zeta
function'', \emph{Dokl. Math.}, \textbf{83}:3 (2011), 314\,–315.

\bibitem{Boyarinov_2011b}

R.N.~Boyarinov, ``On large values of the function $S(t)$ on short
intervals'', \emph{Math. Notes}, \textbf{89}:4 (2011), 472\,–479.

\bibitem{Korolev_2013}

M.A.~Korolev, ``Upper and lower bounds for the function $S(t)$ on
the short intervals'', 2013, 11 pp., arXiv: 1302.0352.

\bibitem{Karatsuba_Korolev_2005}

A.~A.~Karatsuba, M.~A.~Korolev, ``The argument of the Riemann zeta
function'', \emph{Russian Math. Surveys}, \textbf{60}:3 (2005),
433\,-488.

\bibitem{Rosser_Schoenfeld_1962}

J.~B.~Rosser, L.~Schoenfeld, ``Approximate formulas for some
functions of prime numbers'', \emph{Illinois J. Math.}, \textbf{6}:1
(1962), 64\,-94.

\bibitem{Fedoryuk_1977}

M.V.~Fedoryuk, The Saddle-Point Method, Nauka, Moscow, 1977 (in
Russian).

\end{thebibliography}
\end{document}